\documentclass[11pt]{article}
\usepackage{amsmath, amsthm, amssymb, epsfig}
\newtheorem{theorem}{Theorem}[section]
\newtheorem{prop}[theorem]{Proposition}
\newtheorem{lemma}[theorem]{Lemma}
\newtheorem{corollary}[theorem]{Corollary}

\newcommand\beq{\begin{equation}}
\newcommand\eeq{\end{equation}}
\newcommand\bce{\begin{center}}
\newcommand\ece{\end{center}}
\newcommand\bea{\begin{eqnarray}}
\newcommand\eea{\end{eqnarray}}
\newcommand\ben{\begin{enumerate}}
\newcommand\een{\end{enumerate}}
\newcommand\bit{\begin{itemize}}
\newcommand\eit{\end{itemize}}
\newcommand\brr{\begin{array}}
\newcommand\err{\end{array}}
\newcommand\bt{\begin{tabular}}
\newcommand\et{\end{tabular}}

\renewcommand\S{{\mathcal S}}
\newcommand\Al{\operatorname{Allow}}
\newcommand\F{\operatorname{Forb}}
\newcommand\B{\operatorname{B}}
\newcommand\Ba{\operatorname{Bas}}
\newcommand\A{\operatorname{Av}}

\newcommand\Pat{\operatorname{Pat}}
\newcommand\Saw{\operatorname{Saw}}
\newcommand\Poset{P_c}
\newcommand\AS{\operatorname{Alt}}

\title{On basic forbidden patterns of functions}
\author{Sergi Elizalde \hspace{2cm} Yangyang Liu \\ \\
Department of Mathematics\\ Dartmouth College \\ Hanover, NH 03755}
\date{}

\begin{document}
\maketitle
\begin{abstract}
The allowed patterns of a map on a one-dimensional interval are those permutations that are realized by the relative order of the elements in its orbits.
The set of allowed patterns is completely determined by the minimal patterns that are not allowed. These are called basic forbidden patterns.

In this paper we study basic forbidden patterns of several functions. We show that the logistic map $L_r(x)=rx(1-x)$ and some generalizations have infinitely many of them for $1<r\le4$,
and we give a lower bound on the number of basic forbidden patterns of $L_4$ of each length.
Next, we give an upper bound on the length of the shortest forbidden pattern of a piecewise monotone map.
Finally, we provide some necessary conditions for a set of permutations to be the set of basic forbidden patterns of such a map.
\end{abstract}

\section{Introduction and definitions}\label{sec:intro}

Given a map on a one-dimensional interval, consider the finite sequences (orbits) that are obtained by iterating the map,
starting from any point in the interval. The permutations given by the relative order of the elements of these sequences are called {\em allowed patterns};
permutations that do not appear in this way are called {\em forbidden patterns}.
It turns out that piecewise monotone maps always have forbidden patterns, that is,
there are some permutations that do not appear in any orbit~\cite{AEK,BKP}. This idea can be used to distinguish
random sequences, where every permutation appears with some positive probability, from deterministic sequences produced by iterating a map. Practical
aspects of this idea are discussed in~\cite{AZS}.

Minimal forbidden patterns, that is, those for which any proper consecutive subpattern is allowed, are called {\em basic forbidden patterns}. They form an antichain
in the poset of permutations ordered by consecutive pattern containment,
and they contain all the information about the allowed and forbidden patterns of the map.

Consecutive patterns in permutations were first studied in~\cite{EliNoy} from an enumerative point of view. More recently, they have come up in connection to dynamical systems in~\cite{AEK,BKP,Elishifts}.

In this paper we seek to better understand the set of basic forbidden patterns of functions.
Given a map, a natural question is to ask whether its set of basic forbidden patterns is finite or infinite. In Section~\ref{sec:known} we give some easy examples of maps with a finite
set of basic forbidden patterns.
In Section~\ref{sec:logistic} we
show that the set of basic forbidden patterns of the logistic map is infinite, and we find some properties of these patterns.
We show that the result also holds for a more general class of maps.

Section~\ref{sec:shortest} deals with an important practical question. If we are looking for missing patterns in a sequence in order
to tell whether it is random or it has been produced by iterating a piecewise monotone map, it is very useful to have an upper bound on the longest patterns whose presence or absence we need to check.
In Section~\ref{sec:shortest} we provide an upper bound on the length of the shortest forbidden pattern of a map, based on its number of monotonicity intervals.

Another interesting problem is to characterize what sets of permutations can be the basic forbidden patterns of some piecewise monotone map.
In Section~\ref{sec:necessary} we give some necessary conditions that these sets have to satisfy.

\subsection{Permutations and consecutive patterns}
Denote by $\S_{n}$ the set of permutations of $\{1,2,\dots,n\}$. Let $\S=\cup_{n\geq 1} \S_{n}$. If $\pi\in \S_{n}$, we write its one-line notation as $\pi=\pi(1)\pi(2)\dots\pi(n)$. Sometimes it will be convenient to insert commas between the entries.

Let $x_{1},\dots,x_{n}\in \mathbb{R}$ with $x_{1}<x_{2}<\dots<x_{n}$. A permutation of $x_{1},\dots,x_{n}$ can be expressed as $x_{\sigma(1)}x_{\sigma(2)}\dots x_{\sigma(n)}$, where $\sigma\in \S_{n}$.
We define its \emph{reduction} as
$$\rho(x_{\sigma(1)}x_{\sigma(2)}\dots x_{\sigma(n)})=\sigma(1)\sigma(2)\dots\sigma(n)=\sigma.$$ In other words, the reduction is a relabeling of the entries with the numbers $1,2,\dots,n$ while preserving the order relationships among them. For example, $\rho(3,4.2,-2,\sqrt{3},1)=45132$. 

Given two permutations $\pi\in \S_{m}$ and $\sigma\in \S_{n}$ with $m\geq n$, we say that \emph{$\pi$ contains $\sigma$ (as a consecutive pattern)} if there exists $i$ such that $\rho(\pi(i)\pi(i+1)\dots\pi(i+n-1))=\sigma$. In this case, we also say that $\sigma$ is a (consecutive) subpattern of $\pi$, and we write $\sigma\preceq\pi$.
Otherwise, we say that \emph{$\pi$ avoids $\sigma$ (as a consecutive pattern)}.
In the rest of the paper, all the notions of containment, avoidance, and subpatterns refer to the consecutive case.
Denote by $\A_{n}(\sigma)$ the set of permutations in $\S_{n}$ that avoid $\sigma$ as a consecutive pattern, and let $\A(\sigma)=\cup_{n\geq 1}\A_{n}(\sigma)$. In general, if $\Sigma\subset\S$, let $\A(\Sigma)$ be the set of permutations that avoid all the patterns in $\Sigma$, and let $\A_n(\Sigma)=\A(\Sigma)\cap\S_n$. Consecutive pattern containment (and avoidance) was first studied in \cite{EliNoy}.
In~\cite{Eliasym}, the asymptotic behavior of the number of permutations that avoid a consecutive pattern $\pi$ is studied:

\begin{theorem}[\cite{Eliasym}]\label{th:asympconsec}
Let $\sigma\in\S_k$ with $k\ge3$. Then there exist constants $0<c,d<1$ such that $$c^n n! <
|\A_n(\sigma)| < d^n n!$$ for all $n\ge k$.
\end{theorem}

The consecutive containment order $\prec$ defines a partial order on $\S$. Denote by $\Poset$ the corresponding infinite poset.
We say that $A\subset\S$ is a \emph{closed consecutive permutation class} if it is closed under consecutive pattern containment, that is,
if $\pi\in A$ and $\sigma\prec\pi$ imply that $\sigma\in A$.
In this case, the \emph{basis} of $A$ consists of the minimal permutations not in $A$, that is,
$$\Ba(A)=\{\pi\in \S\setminus A: \textrm{if } \sigma\prec\pi, \sigma\neq\pi \textrm{ then } \sigma \in A\}.$$
Note that $\Ba(A)$ is an antichain in $\Poset$ and, conversely, any antichain $\Sigma$ is the basis of the closed class $\A(\Sigma)$.
This gives a one-to-one correspondence between antichains of $\Poset$ and closed consecutive permutation classes.

For example, if $A$ is the set of up-down or down-up permutations, i.e., those permutations satisfying $\pi(1)<\pi(2)>\pi(3)<\pi(4)>\dots$ or $\pi(1)>\pi(2)<\pi(3)>\pi(4)<\dots$,
then $\Ba(A)=\{123,321\}$.
If $B$ is the antichain $\{132,231\}$, then $\A(B)$ is the set
of permutations having no {\em peaks}, i.e., no $i$ such that $\pi(i-1)<\pi(i)>\pi(i+1)$.

\subsection{Allowed and forbidden patterns of maps}
Let $f : I\rightarrow I$, where $I\subset \mathbb{R}$  is a closed interval.
Given $x\in I$ and $n\geq 1$, let $$\Pat(x,f,n)=\rho(x,f(x),f^{2}(x),\dots,f^{n-1}(x)),$$ provided that there is no pair $0\leq i<j\leq n-1$ such that $f^{i}(x)\neq f^{j}(x)$. If such a pair exists, then $\Pat(x,f,n)$ is not defined. When it is defined, $\Pat(x,f,n)\in \S_{n}$. For example, if $L_{4} : [0,1]\rightarrow[0,1]$ is the \emph{logistic map} $L_{4}(x)=4x(1-x)$ and we take $x=0.8$ to be the initial value, then
$$(x,L_4(x),L^2_4(x),L_4^3(x))=(0.8,0.64,0.9216,0.28901376),$$ so $\Pat(0.8,L_4,4)=3241$.

If $\pi\in\S_n$ and there is some $x\in I$ such that $\Pat(x,f,n)=\pi$, we say that $\pi$ is {\em realized} by $f$ (at $x$), or that $\pi$ is an {\em allowed pattern} of $f$.
The set of all permutations realized by $f$ is denoted by $\Al(f)=\bigcup_{n\ge1} \Al_n(f)$, where $$\Al_n(f)=\{\Pat(x,f,n): x\in X\}\subseteq\S_n.$$
The remaining permutations are called {\em forbidden patterns}, and denoted by $\F(f)=\S\setminus\Al(f)$.

It is noticed in~\cite{Elishifts} that $\Al(f)$ is closed under consecutive pattern containment: if $\Pat(x,f,n)=\pi$ and $\tau\prec\pi$, then there exist $i,j$ such that $\rho(\pi(i)\pi(i+1)\dots\pi(j))=\tau$, hence $\Pat(f^{i-1}(x),f,j-i+1)=\tau$.
Those forbidden patterns for which any proper subpattern is allowed are called the {\em basic forbidden patterns} of $f$, and denoted $\B(f)$.
This set is an antichain and it is the basis of $\Al(f)$, i.e., $\B(f)=\Ba(\Al(f))$. In particular, we have that $\Al(f)=\A(\B(f))$.
We will use the notation $\B_{n}(f)=\B(f)\cap\S_n$.

For example, if $g:[0,1]\rightarrow[0,1]$ is the map $g(x)=1-x^2$, then $\B(g)=\{123,132,312,321\}$. To see this, note that the graphs of $x,g(x),g^2(x),\dots$ all intersect at the point
$(\frac{\sqrt{5}+1}{2},\frac{\sqrt{5}+1}{2})$, and that
$$\dots<g^{6}(x)<g^{4}(x)<g^{2}(x)<x<g(x)<g^{3}(x)<g^{5}(x)<\dots$$ for $0<x<\frac{\sqrt{5}+1}{2}$ and
$$\dots>g^{6}(x)>g^{4}(x)>g^{2}(x)>x>g(x)>g^{3}(x)>g^{5}(x)>\dots$$
for $\frac{\sqrt{5}+1}{2}<x<1$.
Other simple cases where the set of basic forbidden patterns is finite and easy to compute are discussed in Section~\ref{sec:known}.

For many maps, however, the set of basic forbidden patterns is infinite.
For the map $L_{4}$ defined above, it can be checked that $\B_{3}(L_{4})=\{321\}$ and $$\B_{4}(L_{4})=\{1423, 2134, 2143, 3142, 4231\}.$$
In Section~\ref{sec:logistic} we study the set of basic forbidden patterns of the logistic map
$$\begin{array}{crcl}
L_{r}: & [0,1]& \longrightarrow & [0,1] \\
&x&\mapsto &r x(1-x),\end{array}$$ where $1<r\leq 4$. We show that $|\B(L_r)|$ is infinite for $1<r\leq 4$, and that $|\B_{n}(L_{4})|\geq n-1$.
We prove that some generalizations of these maps also have infinitely many forbidden patterns.

The following important result, which follows from~\cite{BKP}, guarantees that under mild conditions on $f$, the set $\B(f)$ is nonempty.
Recall that piecewise monotone means that there exists a finite partition of $I$ into intervals where $f$ is continuous and strictly monotone.

\begin{prop}\label{prop1} If  $f : I\rightarrow I$ is piecewise monotone, then $\F(f)\neq\emptyset$. In particular, $\B(f)\neq\emptyset$. \end{prop}

In fact, it is shown in \cite{BKP} that for such a map,
$\lim_{n\rightarrow\infty} \frac{1}{n}\log|\Al_n(f)|$ exists and equals the topological entropy of $f$, a constant which measures the complexity of the dynamical system.
In particular, the number of allowed patterns of $f$ grows at most exponentially, i.e., \beq\label{eq:exponential} |\Al_n(f)|<C^n\eeq for some constant $C$.
Since the total number of permutations of length $n$ grows super-exponentially, the above proposition holds.
In fact, as $n$ approaches infinity, most permutations in $\S_{n}$ are forbidden. In contrast, in a random sequence, all permutations occur with positive probability. Because of this,
forbidden patterns can be used to distinguish random time series from deterministic ones, as studied in \cite{AZS}.

It is shown in~\cite{AEK} that there exist non-piecewise monotone maps that realize all permutations in $\S$. Unless otherwise stated, all the maps $f$ in the rest of the paper will be assumed to be piecewise monotone
maps on an interval $I\subset\mathbb{R}$.

From a practical perspective, the down side of Proposition~\ref{prop1} is that it does not give information about how long the permutations in $\F(f)$ are. Knowing the length of the shortest forbidden pattern of certain classes of maps
is useful when we are trying to distinguish random sequences from chaotic ones generated by orbits of maps in the class. In Section~\ref{sec:shortest} we give an upper bound on the length of the shortest
shortest forbidden pattern of a piecewise monotone map.

Another problem that arises when studying forbidden patterns is the characterization of antichains $\Sigma$ for which there exists a piecewise monotone map $f$ such that $\B(f)=\Sigma$.
This is equivalent to asking whether $\A(\Sigma)$ is the set of allowed patterns of a map, that is,
$\A(\Sigma)=\Al(f)$ for some $f$. It is clear from equation~(\ref{eq:exponential}) that a {\it necessary condition} on $\Sigma$
is that $|\A_n(\Sigma)|<C^n$ for some constant $C$.
For example, if  $\Sigma=\{\sigma\}$, where $\sigma$ has length at least~3, then this condition implies that there is no $f$ such that $\A(\sigma)=\Al(f)$.
Indeed, by Theorem~\ref{th:asympconsec}, $|\A_n(\sigma)|>c^n n!$ for some $0<c<1$, and this lower bound is larger than the necessary exponential growth.
In Section~\ref{sec:necessary} we show that this is not the only necessary condition on the antichain $\Sigma$, that is, there are antichains for
which $|\A_n(\Sigma)|$ grows exponentially, yet there is no piecewise monotone map $f$ with $\B(f)=\Sigma$.

\section{Functions with known forbidden patterns}\label{sec:known}

Determining the forbidden patterns of an arbitrary map is a wide open problem. Only a few results are known for specific maps.
Some work has been done in~\cite{AEK,Elishifts} for the so-called {\em one-sided shift maps}, or simply {\em shifts} for short.
From a forbidden pattern perspective, the shift on $N$ symbols is equivalent to the {\em sawtooth map}
$$\begin{array}{crcl} \Saw_{N}: & [0,1]&\longrightarrow & [0,1] \\
& x & \mapsto & Nx\mod 1,\end{array}$$
as shown in~\cite{AEK}.

It is proved in~\cite{AEK} that shifts (equivalently, sawtooth maps) have infinitely many basic forbidden patterns.
A characterization of forbidden patterns of these maps is given in~\cite{Elishifts},
providing a formula to compute, for a given permutation $\pi$,
the smallest $N$ such that $\pi$ is realized by $\Saw_N$.
The sets $\Al_n(\Saw_N)$ are enumerated for all $n$ and $N$.
To our knowledge, shifts are the only non-trivial maps for which forbidden patterns have been characterized.

\begin{figure}[hbt]
\centering \epsfig{file=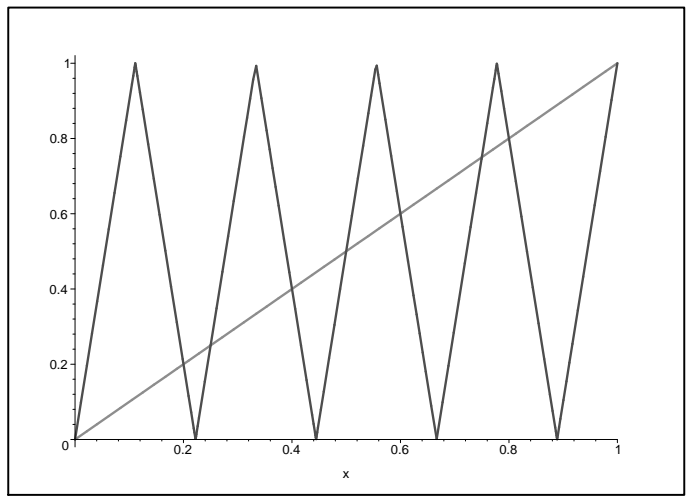,height=9cm,angle=-90} \caption{The alternating sawtooth map $\AS_9$.}\label{fig:altsaw}
\end{figure}

A generalization of shifts are the so-called {\em signed shifts}, which are equivalent to {\em signed sawtooth maps}.
Roughly speaking, for each one of the $N$ spikes of slope $N$ in the graph of $\Saw_N$, one can choose to replace it with a spike
of slope $-N$. For example, for $N=2$, if we reverse the second spike, we obtain the \emph{tent map} $\Lambda:[0,1]\rightarrow[0,1]$ defined by
\begin{equation*} \Lambda (x) = \left\{
\begin{array}{ll}
2x & \text{if } 0\leq x< \frac{1}{2},\\
2-2x & \text{if } \frac{1}{2}\leq x\leq 1.
\end{array} \right.
\end{equation*}
A signed sawtooth map where the slopes alternate between positive and negative is called an {\em alternating signed sawtooth map}. Let $\AS_N:[0,1]\rightarrow[0,1]$ denote the alternating signed sawtooth map with $N$ ramps, defined by
$$\AS_N(x)=\Lambda(\frac{Nx}{2}\mod 1).$$
The graph of $\AS_9$ is shown in Figure~\ref{fig:altsaw}. Forbidden patterns of signed shifts have recently been studied in~\cite{Ami}.

Next we show that for certain well-behaved functions, the description of
their allowed and forbidden patterns is relatively straightforward.

\begin{lemma}\label{lem5} Let $f : [0,1]\rightarrow [0,1]$ be a monotone increasing, continuous function with at least one fixed point on $(0,1)$. Assume that $f$ is not the identity function.
\bit
\item If $f(x)\ge x$ for all $x$, then $\B(f)=\{21\}$;
\item if $f(x)\le x$ for all $x$, then $\B(f)=\{12\}$;
\item otherwise, $\B(f)=\{132,213,231,312\}$.
\eit\end{lemma}

\begin{proof}
Let $U=\{x\in[0,1]:x<f(x)\}$, $V=\{x\in[0,1]:x>f(x)\}$.
Then $U$ decomposes as a union of open intervals $(a,b)$ where $a$ and $b$ are fixed points of $f$, and possibly an interval $[0,b)$.
Since $f$ is increasing and continuous, $f((a,b))=(a,b)$, and $f([0,b))\subseteq[0,b)$. Thus, if $x\in U$ then $f(x)\in U$,
so $x<f(x)<f^{2}(x)<\dots$ for any $x\in U$, and similarly $x>f(x)>f^{2}(x)>\dots$ for any $x\in V$.

If $f(x)\ge x$ for all $x$, then $V=\emptyset$ and the only allowed pattern of length $n$ is $12\dots n$, so $\B(f)=\{21\}$.
Similarly, if $f(x)\le x$ for all $x$, then $U=\emptyset$ and the only allowed pattern of length $n$ is $n(n-1)\dots 1$, so $\B(f)=\{12\}$.
In all other cases, $12\dots n$ and $n(n-1)\dots 1$ are the allowed patterns of length $n$, so $\B(f)=\{132,213,231,312\}$.
\end{proof}

\section{The logistic map and generalizations}\label{sec:logistic}

\subsection{Basic forbidden patterns of the logistic map}

In this section we study the basic forbidden patterns of the logistic map
$$\begin{array}{crcl}
L_{r}(x): & [0,1]& \longrightarrow & [0,1] \\
&x&\mapsto &r x(1-x),\end{array}$$ where $1<r\leq 4$. To simplify notation, we will write $L$ instead of $L_{4}$.

It is shown in~\cite{AEK} that $L$ is order-isomorphic to the tent map $\Lambda$,
and therefore $L$ and $\Lambda$ have the same allowed and forbidden patterns.
It has recently been proved~\cite{Ami} that $\Lambda$ has infinitely many forbidden patterns, by interpreting it as a signed shift.
Here we generalize this result in two ways. First, we show that all maps $L_r$ with $1<r\leq4$ have infinitely many forbidden patterns.
Next, we give a linear lower bound on the number of basic forbidden patterns of $L$ of each length.

\begin{prop}\label{lem2} For every $1<r\leq 4$, $|\B(L_{r})|$ is infinite.
\end{prop}
\begin{proof}
Recall that a permutation of length $n$ is a basic forbidden pattern if and only if it is forbidden and its two subpatterns of length $n-1$ are allowed.

First we show that for $n\ge4$, $(n-2)12\dots(n-3)(n-1)n\in\F(L_{r})$. Let $a_r=1-1/r$ be the unique fixed point of $L_{r}$ in $(0,1)$. It is clear that $\{x\in[0,1]:x< L_{r}(x)\}=(0,a_{r})$ and $\{x\in[0,1]:x>L_{r}(x)\}=(a_{r},1]$.
Suppose, contrary to our claim, that there exists some $x\in[0,1]$ such that $\Pat(x,L_{r},n)=(n-2)123\dots(n-3)(n-1)n$. In other words, $$L_{r}(x)<L_{r}^{2}(x)<\dots<L_{r}^{n-3}(x)<x<L_{r}^{n-2}(x)<L_{r}^{n-1}(x).$$
Then for each $1\leq i\leq n-2$, $L_{r}^{i}(x)\in (0,a_{r})$, whereas $x\in (a_{r},1)$.
However, since $L_{r}^{n-3}(x), L_{r}^{n-2}(x)\in (0,a_{r})$, which is an interval, and $L_{r}^{n-3}(x)<x<L_{r}^{n-2}(x)$, we must have $x\in (0,a_{r})$, leading to contradiction. Therefore, $(n-2)123\dots(n-3)(n-1)n \in \F(L_{r})$.

Next we show that $12\dots m\in\Al(L_r)$ for all $m$. Since $L_{r}(x)=r x(1-x)<4x$ for $x>0$, we have that $L_r^i(x)<a_r/4^{m-2-i}$ for $x\in(0,a_r/4^{m-2})$ and $0\le i\le m-2$.
Now, using that $y<L_r(y)$ for $y\in(0,a_r)$, we get
$$x<L_{r}(x)<L_{r}^{2}(x)<\dots<L_{r}^{m-1}(x)$$
for $x\in(0,a_r/4^{m-2})$, so $123\dots m\in\Al(L_r)$.

Now we show that $m12\dots (m-1)\in\Al(L_r)$ for all $m$. Let $y\in(1-a_r/4^{m-1},1)$ with $y>3/4$. Then $1-y\in(0,{a_r}/{4^{m-1}})$ so, by the same argument as above,
$$1-y<L_{r}(1-y)=L_{r}(y)<L_{r}^{2}(y)<\dots<L_{r}^{m}(y).$$ Also, $L_{r}^{m-1}(y)<a_{r}=\frac{r-1}{r}\le 3/4<y$. Thus $\Pat(y,L_{r},m)=m12\dots(m-1)$.

Summarizing, we have shown that $(n-2)12\dots(n-3)(n-1)n$ is forbidden but both the subpattern formed by its first $n-2$ entries and the one formed by its last $n-1$ entries are allowed.
Now there are two possibilities. If $(n-2)12\dots(n-3)(n-1)$ is forbidden, then it must be a basic forbidden pattern.
If it is allowed, then $(n-2)12\dots(n-3)(n-1)n\in\B(L_r)$. Either way, $\B(L_r)$ contains a permutation of length $n-1$ or $n$, for all $n\ge 4$, so it is an infinite set.
\end{proof}

At the end of the above proof we encountered two possibilities, depending on whether $(n-2)12\dots(n-3)(n-1)$ is forbidden or allowed. We now show that for $r<4$ this pattern is forbidden
for $n$ large enough, but for $r=4$ it is allowed, so $(n-2)12\dots(n-3)(n-1)n\in\B(L)$.

\begin{prop}\label{lem3} For each $1<r<4$, there is some $n_0$ such that for every $n\geq n_0$, $$(n-1)123\dots(n-2)n \in \B(L_{r}).$$ \end{prop}

\begin{proof}
We have shown in the proof of Proposition~\ref{lem2} that $(n-1)12\dots(n-2),123\dots (n-1)\in\Al(L_r)$ for all $n$. So, it suffices to show that $(n-1)123\dots(n-2)n\in\F(L_r)$ for large enough $n$. Let $a_r=1-1/r$ as above.

Let $n\ge3$, and suppose that there exists $x$ such that $\Pat(x,L_{r},n)=(n-1)123\dots(n-2)n$, that is, \beq\label{eqLr}L_{r}(x)<L_{r}^{2}(x)<\dots<L_{r}^{n-2}(x)<x<L_{r}^{n-1}(x).\eeq
Since $x>L_{r}(x)$, we have $x>a_{r}$. Since $L_{r}(x)<L_{r}^{2}(x)$, we have $L_{r}(x)<a_r$, and so $x>\max\{a_r,1-a_r\}=\max\{1-{1}/{r},{1}/{r}\}$.

If $r\leq 2$, then we have $x>{1}/{r}\geq 1/2$. On the other hand, $$\max_{0\leq y\leq 1} L_{r}(y)=\frac{r}{4}\leq \frac{1}{2}.$$ Thus $$x>\max_{0\leq y\leq 1} L_{r}(y)\geq L_{r}^{n-1}(x),$$ contradicting~(\ref{eqLr}).

From now on we assume that $r>2$. In this case $x>1-{1}/{r}=a_r$. Moreover, for $1\le i\le n-3$, we have from~(\ref{eqLr}) that $L_r^{i+1}(x)<L_r^{i+2}(x)$, so $L_r^{i+1}(x)<a_r$, which implies in turn that $L_r^i(x)<1-a_r=1/r$.
Since $L_r(1/r)=1-1/r>1/r$, there is some $\alpha>1$ such that $$L_r(\frac{1}{\alpha r})=\frac{1}{\alpha}(1-\frac{1}{\alpha r})=\frac{1}{r}.$$
Now the fact that
$$L^{n-3}_r(x)<\frac{1}{r}=L_r(\frac{1}{\alpha r})$$
implies that
$L^{n-4}_r(x)<\frac{1}{\alpha r}$ or $L^{n-4}_r(x)>1-\frac{1}{\alpha r}$. Since $L^{n-4}_r(x)<1/r$, it is the first inequality that holds.
By induction on $j$ (we have just done the case $j=1$), we see that
\beq\label{eq:Lj}L^{n-3-j}_r(x)<\frac{1}{\alpha^j r}=\frac{1}{\alpha^{j+1}}(1-\frac{1}{\alpha r})\le \frac{1}{\alpha^{j+1}}(1-\frac{1}{\alpha^{j+1} r})=L_r(\frac{1}{\alpha^{j+1} r}),\eeq
and so, for $j\le n-5$,
$$L^{n-4-j}_r(x)<\frac{1}{\alpha^{j+1} r}.$$
Equation (\ref{eq:Lj}) also holds for $j=n-4$, and from it we get that \beq\label{eq:xalpha} x<\frac{1}{\alpha^{n-3} r} \ \textrm{ or }\ x>1-\frac{1}{\alpha^{n-3} r}.\eeq
Since we know that $x>1-1/r$, the second inequality in~(\ref{eq:xalpha}) must hold.
But since $\alpha>1$ and $r<4$, there must be an $n_0$ such that $$1-\frac{1}{\alpha^{n_0-3} r}>\frac{r}{4}.$$ Now, for $n\ge n_0$,
$$x>1-\frac{1}{\alpha^{n-3} r}>\frac{r}{4}=\max_{0\leq y\leq 1} L_{r}(y)\ge L_{r}^{n-1}(x),$$ again contradicting~(\ref{eqLr}).
\end{proof}

For $r=4$, forbidden patterns behave in a different way. The pattern mentioned in Proposition~\ref{lem3} is now allowed, but instead we have other basic forbidden patterns.
It is shown in~\cite[Theorem~4.4]{Ami} that $(n-3)(n-2)(n-1)12\dots(n-4)n\in\B_n(L)$ for $n\ge5$.
We find here $n-1$ additional basic forbidden patterns of length $n$, thus giving a linear lower bound on $|\B_n(L)|$.

\begin{prop}\label{lem4} For $n\geq 4$, the set $\B_{n}(L)$ contains the following patterns:
$$\begin{array}{l} (n-2)12\dots(n-3)(n-1)n,\\
(n-2)12\dots(n-3)n(n-1),\\
(n-1)12\dots(k-1)(k+1)\dots(n-2)nk \quad \textrm{ for }2\leq k\leq n-2.\end{array}$$
In particular
 $$|\B_{n}(L)|\geq n.$$
\end{prop}

\begin{proof}
It is shown in~\cite{AEK} that $\B(L)=\B(\Lambda)$, so it suffices to prove the statement for the tent map $\Lambda$.

First we show that $(m-1)123\dots(m-2)m\in\Al(\Lambda)$.
Let $$x=1-\frac{1+2^{-m}}{2^{m-1}+1}.$$
Then $\Lambda(x)=2(1-x)$. In general, for $1\le i\le m-2$,
$$\Lambda^{i}(x)=2^i(1-x)=\frac{2^i+2^{i-m}}{2^{m-1}+1}\le \frac{2^{m-2}+1/4}{2^{m-1}+1}<\frac{1}{2}<x,$$
and $$\Lambda^{m-1}(x)=2^{m-1}(1-x)=\frac{2^{m-1}+1/2}{2^{m-1}+1}>x.$$ Thus,
$$\Lambda(x)<\Lambda^2(x)<\dots<\Lambda^{m-2}(x)<x<\Lambda^{m-1}(x),$$
so $\Pat(x,\Lambda,m)= (m-1)123\dots(m-2)m$.

Now we show that $123\dots(k-1)(k+1)\dots(m-1)mk\in \Al(\Lambda)$ for $1\leq k\leq m$.
Since $$\{x\in[0,1]:x<\Lambda(x)\}=(0,{2}/{3})$$ and $$\{x\in[0,1]:\Lambda(x)<c\}=(0,c/2)\cup(1-c/2,1)$$ for any $0<c<1$, we have that
$$\{x\in [0,1]:x<\Lambda(x)<\Lambda^{2}(x)\}=(0,{1}/{3}).$$
In general, it is easy to see that for $m\ge3$, \beq\label{eq:incrlambda}\{x\in[0,1]:x<\Lambda(x)<\Lambda^{2}(x)<\dots<\Lambda^{m-2}(x)\}=(0,\frac{1}{3\cdot 2^{m-4}}).\eeq
Since $\Lambda^{m-1}(x)$ is continuous, $\Lambda^{m-1}(1/2^{m-1})=1$ and $\Lambda^{m-1}({1}/{2^{m-2}})=0$, the graph of $\Lambda^{m-1}(x)$ intersects
the graph of $\Lambda^{i}(x)$, for each $0\leq i\leq m-2$, in the interval $$\left[\frac{1}{2^{m-1}},\frac{1}{2^{m-2}}\right]\subset\left(0,\frac{1}{3\cdot 2^{m-4}}\right),$$ thus realizing the patterns $123\dots(k-1)(k+1)\dots mk$ for each $1\leq k\leq m$.

Since for each of the patterns in the statement of the proposition both of its subpatterns of length $n-1$ are allowed, it suffices to show that they are forbidden to conclude that they are in $\B(f)$.

Suppose there exists $x\in[0,1]$ such that $\Pat(x,\Lambda,n)$ is one of the listed patterns. Then $$\Lambda(x)<\Lambda^{2}(x)<\dots<\Lambda^{n-3}(x)<x<\Lambda^{n-2}(x).$$
Since $x>\Lambda(x)$, we have $x>2/3$. And since $\Lambda^{n-2}(x)>x>1/2$, we have that $$\Lambda^{n-1}(x)=2(1-\Lambda^{n-2}(x))<2(1-x)=\Lambda(x).$$
Therefore, $$\Lambda^{n-1}(x)<\Lambda(x)<\Lambda^{2}(x)<\dots<\Lambda^{n-3}(x)<x<\Lambda^{n-2}(x),$$ so
$\Pat(x,\Lambda,n)=(n-1)23\dots(n-2)n1$, and all the patterns in the statement are forbidden.

Together with the fact that $(n-3)(n-2)(n-1)12\dots(n-4)n\in\B_n(L)$ for $n\ge5$ and that $|\B_{4}(L)|=5$, the lower bound $|\B_n(L)|\ge n$ follows.
\end{proof}

In fact we expect the actual size of $\B_n(L)$ to grow much faster than this. For $n\ge3$, the first few
values of $|\B_n(L)|$, found by computer, are $1,5,9,28,53,110,\dots$.

\subsection{Generalizations}
It is an interesting open problem to characterize those maps for which the set of basic forbidden patterns in infinite. We showed that this is the case for $L_r$. Here we give
a sufficient set of conditions on $f$ that makes $\B(f)$ infinite. The conditions generalize some properties of the logistic map that we used in the above subsection, including symmetry
and having a single fixed point in $(0,1)$.

\begin{prop}\label{prop5} Let $f:[0,1]\rightarrow[0,1]$ be a continuous function that satisfies the following three conditions:
\ben 
\item $f(0)=0$,
\item $f(x)=f(1-x)$ for all $x$,
\item $f(x)$ has a single fixed point in $(0,1)$,
\item there is some $x$ such that $f(x)>x$.
\een
Then $f$ has infinitely many basic forbidden patterns.
\end{prop}

Note that in the above proposition we do not require $f$ to be piecewise monotone.

\begin{proof}
This proof is similar to that of Proposition~\ref{lem2}.
For each $i\geq 1$, let $$T_{i}= \{x\in[0,1]:x<f(x)<f^{2}(x)<\dots<f^{i}(x)\}.$$
Conditions 1, 3 and 4 imply that $T_{1}=(0,t_{1})$, where $t_{1}$ is the fixed point of $f$ in $(0,1)$.

We claim that there exists a decreasing sequence of positive numbers $t_{1}>t_{2}>t_{3}>\dots$ such that $(0,t_{i})\subseteq T_{i}$ for each $i\geq 1$.
We show that $t_i$ exists by induction on $i$. Assume that $i\ge2$ and that we have shown the existence of $t_{i-1}$ with $(0,t_{i-1})\subseteq T_{i-1}$.
Let $t_{i}$ be the smallest root of $f(x)=t_{i-1}$. Clearly $t_{i}\leq t_{i-1}$. Let $x\in (0,t_{i})$. Since $f(0)=0<t_{i-1}$, we have that $f(x)<t_{i-1}$, so $f(x)\in T_{i-1}$. By the definition of $T_{i-1}$,
$f(x)<f^{2}(x)<\dots<f^{i}(x)$. Moreover, since $x<t_{i}\le t_{i-1}$, $x\in T_{i-1}$, so $x<f(x)$. It follows that $x\in T_{i}$. Hence $(0,t_{i})\subseteq T_{i}$.
This proves that for all $m$, $T_{m-1}\neq\emptyset$, so $12\dots m\in\Al(f)$.

We next show that for $m\geq 3$, $m123\dots(m-1)\in\Al(f)$. Let $x\in(\max\{t_{1}, 1-t_{m}\},1)$. Then $1-x\in T_{m}$, so $$1-x<f(1-x)=f(x)<f^2(x)<\dots<f^m(x).$$
Since $f^{m-1}(x)<f^{m}(x)$, we have that $f^{m-1}(x)\in T_{1}$, so $f^{m-1}(x)<t_{1}<x$. Thus $$f(x)<f^{2}(x)<\dots<f^{m-1}(x)<x,$$ which shows that $m123\dots(m-1)\in\Al(f)$.

Now we prove that for $n\ge4$, $(n-2)12\dots(n-3)(n-1)n\in \F(f)$. Suppose that there exists some $x\in [0,1]$ such that $\Pat(x,f,n)=(n-2)12\dots(n-3)(n-1)n$. Then $x>f(x)$ and $f^{n-3}(x)<f^{n-2}(x)<f^{n-1}(x)$.
Therefore, $x\notin T_{1}$ whereas $f^{n-3}(x),f^{n-2}(x)\in T_{1}$. However, because $T_{1}=(0,t_{1})$ is an interval and $f^{n-3}(x)<x<f^{n-2}(x)$, we must have $x\in T_{1}$, which is a contradiction.

Summarizing, we have shown that $(n-2)12\dots(n-3)(n-1)n$ is forbidden but both the subpattern formed by its first $n-2$ entries and the one formed by its last $n-1$ entries are allowed.
If $(n-2)12\dots(n-3)(n-1)$ is forbidden, then it must be a basic forbidden pattern; if it is allowed, then $(n-2)12\dots(n-3)(n-1)n\in\B(f)$.
Either way, $\B(f)$ contains a permutation of length $n-1$ or $n$, for all $n\ge 4$.
\end{proof}

Obviously, the conditions in Proposition~\ref{prop5} are not necessary for a map to have infinitely many forbidden patterns. For example, as mentioned above, it is known~\cite{AEK} that the maps $\Saw_N$ for $N\ge2$ have infinitely many forbidden patterns.

\section{The shortest forbidden pattern}\label{sec:shortest}

In this section we give an upper bound on the length of the shortest forbidden pattern of a piecewise monotone map $f:I\rightarrow I$.
The intervals in the partition of $I$ such that $f$ is continuous and strictly monotone on each interval are called the {\em monotonicity intervals} of~$f$.

\begin{theorem}\label{thm:shortest} Let $f:I\rightarrow I$ be a piecewise monotone map with $m$ monotonicity intervals, and let $k\le m$ be the number of such intervals $I'$ with one of these two properties:
\bit \item $f$ is increasing in $I'$ and the left endpoint $a$ of $\overline{I'}$ satisfies $f(a)<a$,
\item $f$ is decreasing in $I'$ and $I'$ contains a point with $f(x)<x$.
\eit
Then the length of the shortest forbidden pattern of $f$ is at most $2k+3$.\end{theorem}

\begin{proof}
Let $I=\bigcup_{j=1}^m I_j$ be the finite partition into intervals where $f$ is continuous and strictly monotone.
Let $$D=\{x\in I:f(x)<x\}=\bigcup_\alpha D_\alpha$$ expressed as the union of its connected components,
where each $D_\alpha$ is an interval (which can consist of a single point).
Note that this can be an infinite union.

Let $1\le j\le m$. If $f$ is decreasing in $I_j$, then $I_j$ can intersect at most one of the $D_\alpha$.
Suppose now that $f$ is increasing in $I_j$, and let $\overline{I_j}=[a_j,b_j]$ and $\overline{D_\alpha}=[c_\alpha,d_\alpha]$.
We claim that if $a_j<c_\alpha<b_j$, then $f(D_\alpha)\subseteq D_\alpha$.
Indeed, since $f$ is continuous in $I_j$ and $c_\alpha$ is the left endpoint of $D_\alpha$, we have that $f(c_\alpha)=c_\alpha$,
and so $$c_\alpha<f(x)<x\le d_\alpha$$ for all $x\in D_\alpha$, which implies that $f(D_\alpha)\subseteq D_\alpha$.
For the same reason, if $a_j=c_\alpha$ but $f(a_j)=a_j$, we have again that $f(D_\alpha)\subseteq D_\alpha$.
Therefore, the number of intervals $D_\alpha$ for which $f(D_\alpha)\nsubseteq D_\alpha$ is at most $k$.

For $n\ge 2k+3$, let $E_n\subset\S_n$ be the set of permutations $\pi$ for which there exists an $i$ such that 
\ben
\item for each $i\le j\le i+k-1$, there exists $1\le\ell\le n-1$ such that $\pi(j)>\pi(\ell)>\pi(j+1)$ and $\pi(\ell)<\pi(\ell+1)$ (note that this condition implies that $\pi(j)\ge\pi(j+1)+2$),
\item $\pi(i+k)>\pi(i+k+1)$,
\item there is some $h>i+k+1$ such that $\pi(i+k+1)<\pi(h)$.
\een
For example, $\pi=(2k+2)(2k)(2k-2)\dots42135\dots(2k+1)(2k+3)\in E_{2k+3}$.
We claim that $E_n\subset \F(f)$. Once we prove this claim, the theorem will follow from this example.

Suppose that for some $\sigma\in E_n$ there is a $y\in I$ such that $\Pat(y,f,n)=\sigma$. Then,
whenever $\pi(j)>\pi(j+1)$, we have $f^{j-1}(y)>f^j(y)$ and thus $f^{j-1}(y)\in D$, and
whenever $\pi(\ell)<\pi(\ell+1)$, we have $f^{\ell-1}(y)<f^{\ell}(y)$ and thus $f^{\ell-1}(y)\notin D$.

For each $i\le j\le i+k$, we have $f^{j-1}(y)\in D_{\alpha_j}$ for some $\alpha_j$. If $f(D_{\alpha_j})\subseteq D_{\alpha_j}$, then $f^{j-1}(y)>f^j(y)>f^{j+1}(y)>\dots$. But this
is impossible because condition 3 implies that $f^{i+k}(y)<f^{h-1}(y)$ for some $h>i+k+1$. Therefore, $f(D_{\alpha_j})\nsubseteq D_{\alpha_j}$.
Now, there are $k+1$ choices for $j$, and only $k$ different indices $\alpha$ such that $f(D_\alpha)\nsubseteq D_\alpha$. Hence,
there must be two different indices $i\le j<j'\le i+k$ for which $\alpha_j=\alpha_{j'}$.
Since $D_{\alpha_j}$ is an interval, for any $f^{j-1}(y)>z>f^{j'-1}(y)$ we must have $z\in D_{\alpha_j}\subseteq D$.
However, by condition 1, there exists $\ell$ such that $\pi(j)>\pi(\ell)>\pi(j+1)\ge\pi(j')$ and $\pi(\ell)<\pi(\ell+1)$.
This implies that $f^{j-1}(y)>f^{\ell-1}(y)>f^{j'-1}(y)$ but at the same time $f^{\ell-1}(y)<f^{\ell}(y)$, so $f^{\ell-1}(y)\notin D$, which is a contradiction.
\end{proof}

The permutation $\pi$ in the above proof is not the only element of $E_{2k+3}$. For example, for $k=1$, we have $34215,35214,42135,45213,45312,52134\in E_5$.

In many cases, it is more practical to work with the following simplified version of Theorem~\ref{thm:shortest}.
\begin{theorem}\label{thm:shortesteasy} Let $f:I\rightarrow I$ be a piecewise monotone map and let $D=\{x\in I: f(x)<x\}$. Let $k$ be the number of connected components of $D$.
Then the length of the shortest forbidden pattern of $f$ is at most $2k+2$.\end{theorem}

\begin{proof}
This proof is analogous to that of Theorem~\ref{thm:shortest}. In this case, the statement already gives the bound $k$ on the number of intervals $D_\alpha$. Since now we do not need to eliminate those with
$f(D_\alpha)\subseteq D_\alpha$, condition 3 in the definition of $E_n$ can be dropped. We have that $E'_n\subset \F(f)$ for $n\ge 2k+2$, where
$E'_n\subset\S_n$ be the set of permutations $\pi$ for which there exists an $i$ such that
\ben
\item for each $i\le j\le i+k-1$, there exists $1\le\ell\le n-1$ such that $\pi(j)>\pi(\ell)>\pi(j+1)$ and $\pi(\ell)<\pi(\ell+1)$,
\item $\pi(i+k)>\pi(i+k+1)$.
\een
For example, $\pi=35\dots(2k+1)(2k+2)(2k)(2k-2)\dots421\in E'_{2k+2}$.
\end{proof}

We can apply Theorem~\ref{thm:shortesteasy} to the sawtooth map $\Saw_N$. In this case, $D$ is the union of $k=N-1$ intervals, and the theorem guarantees that the
shortest forbidden pattern has length at most $2N$. In fact, it is shown in~\cite{AEK} that the shortest forbidden pattern of $\Saw_{N}$ has length $N+2$.
For small values of $N$, it follows from the proof of Theorem~\ref{thm:shortesteasy} that $3421\in\F(\Saw_2)$ and $356421\in\F(\Saw_3)$.

Our theorem gives a tight bound when applied to some alternating sawtooth maps. For the map $\AS_N$ where $N$ is odd, we see that $D$ has $k=(N-1)/2$ components (see Figure~\ref{fig:altsaw}).
In this case, Theorem~\ref{thm:shortesteasy} states that the shortest forbidden pattern of $\AS_N$ has length at most $N+1$, and this turns out to be its actual length, as shown in~\cite[Theorem 4.5]{Ami}.

Both Theorem~\ref{thm:shortest} and~\ref{thm:shortesteasy} have analogous symmetric formulations if we consider the set $\{x\in I: f(x)>x\}$ instead of $D$.
For example, here is the corresponding version of Theorem~\ref{thm:shortest}.

\begin{corollary} Let $f$ be a piecewise monotone map with $m$ monotonicity intervals. Let $k\le m$ be the number of such intervals $I'$ with one of these two properties:
\bit \item $f$ is increasing in $I'$ and the right endpoint $a$ of $I'$ satisfies $f(a)>a$,
\item $f$ is decreasing in $I'$ and $I'$ contains a point with $f(x)>x$.
\eit
Then the length of the shortest forbidden pattern of $f$ is at most $2k+3$.\end{corollary}

\section{Antichains that are basic forbidden patterns of a function}\label{sec:necessary}

In Section~\ref{sec:intro} we mentioned the problem of characterizing those sets $\Sigma$ for which there exists a piecewise monotone map $f$ such that $\B(f)=\Sigma$.
Aside from the obvious prerequisite that $\Sigma$ has to be an antichain in $\Poset$, another necessary condition is that the number of permutations avoiding $\Sigma$ must grow at most exponentially.
Using this requirement we can show that certain finite antichains are not of the form $\B(f)$. In the next proposition, the floor function $[x]$ is the largest integer that is less than or equal to $x$.

\begin{prop}\label{prop3} Let $\Sigma=\{\sigma_{1}, \sigma_{2}, \dots, \sigma_{m}\}$ be a finite antichain in $\Poset$ and let $k_{i}$ be the length of $\sigma_{i}$, with $k_{1}\leq k_{2} \leq \dots \leq k_{m}$.
If $B$ is the set of basic forbidden patterns of a piecewise monotone map, then $$k_{1}+k_{2}+\dots+k_{m}\geq [k_{1}/2]!+m([k_{1}/2]-1).$$
\end{prop}

\begin{proof}
Let $\ell=[k_{1}/2]$. Assume to the contrary that $k_{1}+k_{2}+\dots+k_{m}< \ell!+m(\ell-1)$. Equivalently, $$\sum_{i=1}^{m} (k_{i}-\ell+1) < \ell!.$$
There are $k_{i}-\ell+1$ consecutive subpatterns (not necessarily different) of length $\ell$ in $\sigma_{i}$. Thus,
the above inequality implies that there is at least one permutation $\pi$ of length $\ell$ that is not contained in any of $\sigma_{1}, \sigma_{2}, \dots, \sigma_{m}$.

Let
$$G_{r\ell}=\{\tau_{1}\tau_{2}\dots\tau_{r}\in\S_{r\ell}: \rho(\tau_i)=\pi \textrm{ for } 1\le i\le r\}.$$  We claim that $G_{r\ell}\subseteq \A(\Sigma)$. To see this, note that every subpattern of $\tau_{1}\tau_{2}\dots\tau_{r}$
of length at least $k_{1}$ spans the entirety of some $\tau_{i}$, so it contains $\pi$. On the other hand, no permutation in $\Sigma$ contains $\pi$.
Therefore, no permutation in $\Sigma$ is a subpattern of any $\tau_{1}\tau_{2}\dots\tau_{r}\in G_{r\ell}$, so $\tau_{1}\tau_{2}\dots\tau_{r}$ avoids $\Sigma$.

The size of $G_{r\ell}$ is equal to the number of ways to partition the set $\{1,2,\dots,r\ell\}$ into $r$ blocks of size $\ell$, which is
$$|G_{r\ell}|=\binom{r\ell}{\ell,\ell,\dots,\ell}=\frac{(r\ell)!}{(\ell!)^{r}}.$$
Using Stirling's formula, we see that as $r$ goes to infinity,
$$|\A_{r\ell}(\Sigma)|\ge|G_{r\ell}|\gg \frac{(r\ell)^{r\ell}}{e^{r\ell}{(\ell!)^{r}}}=\left(\frac{\ell^{\ell}}{e^{\ell}\ell!}\right)^{r}r^r,$$
so $|\A_{r\ell}(\Sigma)|$ grows super-exponentially.
Thus $\Sigma$ cannot be the set of basic forbidden patterns of piecewise monotone map.
\end{proof}

Theorem~\ref{th:asympconsec} states that if $\sigma$ has length $k\geq 3$, $|\A_{n}(\sigma)|$ grows super-exponentially. When $k\geq 6$, this result can be directly derived from the above proposition. To see this, note that when $k_{1}\geq 6, [k_{1}/2]\geq 3$, and so
$$[k_{1}/2]!+n([k_{1}/2]-1)-\sum_{i=1}^{n} k_{i}=[k_{1}/2]!+[k_{1}/2]-1-k_{1}
\geq 2[k_{1}/2]+[k_{1}/2]-1-k_{1}>0.$$

Interestingly, exponential growth on the number of permutations avoiding an antichain is not the only requirement for it to be the set of basic forbidden patterns of a piecewise monotone function.
For example, consider the antichain $\Sigma=\{132, 231\}$. Then, as mentioned in the introduction, $\A(\Sigma)$ is the set of permutations with no peaks.
In other words, permutations in $\A(\Sigma)$ consist of a decreasing sequence followed by an increasing sequence. It is easy to see that $|\A_n(\Sigma)|=2^{n-1}$, since such a permutation is
determined by the set of elements other than $1$ in the initial decreasing sequence. However, even though the exponential growth condition is satisfied, we have the following result.

\begin{prop}\label{lem1} There exists no piecewise monotone map $f$ on a closed interval $I\subset\mathbb{R}$ such that $\B(f)=\{132, 231\}$. \end{prop}

\begin{proof}
Assume to the contrary that there exists such a map $f$. Let $m$ be the number of monotonicity intervals of $f$. As in the proof of Theorem~\ref{thm:shortest},
we have that $\pi=(2m+2)(2m)(2m-2)\dots42135\dots(2m-1)(2m+1)(2m+3)\in \F(f)$.
Since $\pi$ avoids $132$ and $231$, it is not possible that $\B(f)=\{132, 231\}$.
\end{proof}


\begin{thebibliography}{}
\bibitem{Ami} J.M. Amig\'o, The ordinal structure of the signed shift transformations, {\it International Journal of Bifurcation and Chaos}, to appear.

\bibitem{AEK} J.M. Amig\'{o}, S. Elizalde, and M. Kennel, Forbidden patterns and shift systems, \textit{J. Combin. Theory Ser. A} 115 (2008), 485--504.

\bibitem{AZS} J.M. Amig\'{o}, S. Zambrano, and M.A.F. Sanju\'{a}n, True and false forbidden patterns in deterministic and random dynamics, \textit{Europhys. Lett.} 79 (2007), 50001-p1, -p5.

\bibitem{BKP} C. Bandt, G. Keller, and B. Pompe, Entropy of interval maps via permutations, \textit{Nonlinearity} 15 (2002), 1595--1602.

\bibitem{Eliasym} S. Elizalde, Asymptotic enumeration of permutations
avoiding generalized patterns, {\it Adv. in Appl. Math.}
36 (2006), 138--155.

\bibitem{Elishifts} S. Elizalde, The number of permutations realized by a shift, \textit{SIAM J. Discrete Math.} 23 (2009), 765--786.

\bibitem{EliNoy} S. Elizalde and M. Noy, Consecutive patterns in permutations, \textit{Adv. Appl. Math.} 30 (2003), 110--123.
\end{thebibliography}
\end{document}